\documentclass{amsart}
\usepackage{graphicx,latexsym,amssymb}

\newtheorem{thm}{Theorem}[section]
\newtheorem{lmm}{Lemma}[section]
\newtheorem{cor}{Corollary}[section]
\newtheorem{prop}[thm]{Proposition}
\newtheorem{defn}[thm]{Definition}
\newtheorem{rem}[thm]{Remark}
\newtheorem{example}[thm]{Example}

\newcommand{\C}{\mathbb{C}}

\newcommand{\CP}{\mathbb{CP}}

\newcommand{\ben}{\begin{enumerate}}
\newcommand{\een}{\end{enumerate}}
\newcommand{\ble}{\begin{lem}}
\newcommand{\ele}{\end{lem}}
\newcommand{\bth}{\begin{thm}}
\renewcommand{\eth}{\end{thm}}
\newcommand{\bpr}{\begin{prop}}
\newcommand{\epr}{\end{prop}}
\newcommand{\bco}{\begin{cor}}
\newcommand{\eco}{\end{cor}}
\newcommand{\bde}{\begin{defn}}
\newcommand{\ede}{\end{defn}}
\newcommand{\brem}{\begin{rem}}
\newcommand{\erem}{\end{rem}}
\newcommand{\bexm}{\begin{example}}
\newcommand{\eexm}{\end{example}}

\begin{document}

\title{Rational transformations of algebraic curves and elimination theory}

\author{Alexander Shapiro and Victor Vinnikov}

\address{Alexander Shapiro,
Department of Mathematics and Statistics, Bar-Ilan University, Ramat-Gan 52900, Israel}
\email{sapial@math.biu.ac.il}

\address{Victor Vinnikov,
Department of Mathematics, Ben-Gurion University of the Negev, Beer-Sheva 84105, Israel}
\email{vinnikov@cs.bgu.ac.il}

\thanks{Partially supported by  EU-network HPRN-CT-2009-00099(EAGER) , (The Emmy Noether Research Institute for Mathematics and the Minerva Foundation of Germany), the Israel Science Foundation grant \# 8008/02-3 (Excellency Center "Group  Theoretic Methods in the Study of Algebraic Varieties").}

\maketitle

\begin{abstract}
Elimination theory has many applications, in particular, it describes explicitly an image of a complex line
under rational transformation and determines the number of common zeroes of two polynomials in one
variable. We generalize classical elimination theory and create elimination theory along an algebraic curve
using the notion of determinantal representation of algebraic curve. This new theory allows to describe
explicitly an image of a plane algebraic curve under rational transformation and to determine the number of
common zeroes of two polynomials in two variables on a plane algebraic curve.
\end{abstract}

\section*{Introduction}

The main goal of this research is to describe explicitly an image of an algebraic curve under rational
transformation.

The simplest and very illustrative case is a rational transformation of a projective line into projective
plane. Three homogeneous polynomials in two variables $p_0$, $p_1$ and $p_2$ maps a projective line $\CP$
into projective plane $\CP^2$:
$$(x_0,x_1) \to (p_0(x_0,x_1), p_1(x_0,x_1), p_2(x_0,x_1))$$

This case was described by N.~Kravitsky using the classical elimination theory, see \cite{K}. The image of
a projective line is a rational curve. This curve is defined by a polynomial
$${\Delta}(x_0,x_1,x_2)=
\det(x_0B(p_1,p_2)+x_1B(p_2,p_0)+x_2B(p_0,p_1))$$
where $B(p_i,p_j)$ is the Bezout matrix of polynomials $p_i$ and $p_j$.

Our original objective was to find an analogue of the constructions of \cite{K} in the general case. This
led us to consider elimination theory for pairs of polynomials along an algebraic curve given by a
determinantal representation. While our results, as presented in this paper, are for polynomials in two
variables, plane algebraic curves, and pairs of operators, the generalization to polynomials in $d$
variables, algebraic curves in the $d$-dimensional space, and $d$-tuples of operators should be, for the
most part, relatively straightforward.

Let us recall the main goal of elimination theory. Given $n+1$ (nonhomogeneous) polynomials in $n$
variables we want to find necessary and sufficient conditions (in terms of the coefficients) for these
polynomials to have a common zero (and furthermore to determine the number of common zeroes, counting
multiplicities, if they exist), see \cite{M}.

In the classical case we consider (nonhomogeneous) polynomials in one
variable,
$$p(x) = p_0 + p_1 x + p_2 x^2+ \dots + p_n x^n =
\sum_{i=0}^n p_i x^i =$$
$$=(p_0, p_1, \dots , p_n)
\cdot
\left(
\begin{array}{c}
1 \\
x \\
x^2 \\
\dots \\
x^n
\end{array}
\right)=pV_{n+1}(x)$$
Here $p=(p_0, p_1, \dots , p_n)$ is the row--vector of coefficients of
the polynomial $p(x)$ and $V_{n+1}(x)$ is a so called Vandermonde
vector of the length $n+1$:
$$V_{n+1}(x)=
\left(
\begin{array}{c}
1 \\
x \\
x^2 \\
\dots \\
x^n
\end{array}
\right)$$ We will see that Vandermonde vectors provide a convenient framework for a proof of classical
results of elimination theory for two polynomials in one variable. In fact a key point in generalizing the
classical results in our approach is to find an appropriate generalization of the notion a Vandermonde
vector.

For polynomials in two variables,
$$p(x_1,x_2) = p_{00}+p_{10}x_1+p_{01}x_2+\dots
+p_{n0}x_1^n+p_{n-11}x_1^{n-1}x_2+ \dots +p_{0n}x_2^n=$$
$$=\sum_{i_1,i_2=0}^n p_{i_1,i_2}x_1^{i_1}x_2^{i_2}=$$
$$=\left(\begin{array}{c}
p_{00}\\
p_{10}{\quad}p_{01}\\
p_{20}{\quad}p_{11}{\quad}p_{02}\\
\dots{\quad}\dots{\quad}\dots{\quad}\dots\\
p_{n0}{\quad}p_{n-11}{\quad}\dots{\quad}\dots p_{0n}
\end{array}\right)^T
\cdot
\left(\begin{array}{c}
1\\
x_1{\quad}x_2\\
x_1^2{\quad}x_1x_2{\quad}x_2^2\\
\dots{\quad}\dots{\quad}\dots{\quad}\dots\\
x_1^n{\quad}x_1^{n-1}x_2{\quad}\dots{\quad}\dots{\quad}x_2^n
\end{array}\right),$$ one usually considers elimination theory
for three polynomials.

But it is an unfortunate fact that the easy, useful and beautiful constructions of the single variable case
do not generalize in any straightforward fashion to the case of two or more variables. Some generalizations
of the classical results using much more sophisticated algebraic techniques have been obtained in recent
years by Gelfand, Kapranov, Zelevinsky, see \cite{GKZ}, and Jouanolou, see \cite{J}. We shall follow a
different, somewhat ``asymmetrical'' approach. Namely, we shall choose one of our three polynomials and
view it as defining a plane algebraic curve.

Our goal becomes to determine, in terms of coefficients of two polynomials in two variables, whether these
two polynomials have a common zero on an algebraic curve defined by a third polynomial, and if they do to
determine the number of common zeroes, counting multiplicities.

To solve this problem it turns out to be essential to consider a plane algebraic curve defined by a
polynomial $\Delta(x_1,x_2)$ of degree $m$ together with a determinantal representation of the form
$${\Delta}(x_1,x_2)=\det(D_0+x_1D_1+x_2D_2)$$
where $D_0$, $D_1$ and $D_2$ are $m \times m$ constant matrices. In this paper we shall assume that the
polynomial ${\Delta}(x_1,x_2)$ is irreducible, though many results can be generalized to the reducible case
as well. We shall also assume that ${\Delta}(x_1,x_2)$ has real coefficients and $D_0$, $D_1$ and $D_2$ are
complex hermitian matrices; it allows to dispose of certain considerations of duality which will be
necessary to handle the case when the reality assumption is dropped.

The determinantal representation of a curve provides us with an additional structure. In the case of a
plane irreducible curve for every point $(x_1,x_2)$ of a curve there is a nontrivial subspace
$\ker(D_0+x_1D_1+x_2D_2)$ (one-dimensional, because of the irreducibility of $\Delta(x_1,x_2)$, except at
the singular points) -- the fiber of the corresponding line bundle at this point. We choose a vector $e$ in
this fiber. This additional structure allows us to define an appropriate analogue of Vandermonde vectors in
one variable: Vandermonde vectors on a curve. Vandermonde vector on a curve (equipped with a determinantal
representation as above) is given by
$$V_{n+1}(x_1,x_2,e)=
\left(\begin{array}{c}
e\\
x_1e{\quad}x_2e\\
x_1^2e{\quad}x_1x_2e{\quad}x_2^2e\\
\dots{\quad}\dots{\quad}\dots{\quad}\dots\\
x_1^ne{\quad}x_1^{n-1}x_2e{\quad}\dots{\quad}\dots{\quad}x_2^ne
\end{array}\right)$$

Let $m$ be the degree of the polynomial ${\Delta}(x_1,x_2)=\det(D_0+x_1D_1+x_2D_2)$. We will denote the
complex euclidian space of the dimension $m\frac{n(n+1)}{2}$ by $W_n$. It is clear that every Vandermonde
vector $V_n(x_1,x_2,e)$ on a curve belongs to $W_n$. We will call $W_n$ the blown space. Unlike the one
variable case all the Vandermonde vectors on a curve generate not the whole blown space $W_n$ but only a
certain subspace of $W_n$. Let us denote the subspace generated by Vandermonde vectors on a curve by $V_n$.
We will call $V_n$ the principal subspace. The dimension of the principal subspace $V_n$ is equal to $nm$,
and it can be described explicitly in terms of the matrices $D_0$, $D_1$ and $D_2$ appearing in the
determinantal representation.

The blown space and the principal subspace play a key role in the elimination theory on a plane algebraic
curve. Namely, all the constructions of the classical elimination theory can be formally generalized in
terms of the blown spaces. After this formal generalization we restrict the new ``blown'' constructions to
the corresponding principal subspaces and obtain complete analogues of the classical results.

One immediate result of elimination theory along a plane algebraic curve is that these analogues of the
classical constructions allow us to describe an image of an algebraic curve given by a determinantal
representation under a rational transformation.

As we have mentioned earlier we expect that all the results of our elimination theory extend in a
straightforward fashion to algebraic curves in a $d$-dimensional space for any $d$, using determinantal
representations of the form discussed in \cite{LMKV}, p. 42, yielding a properly defined notion of a
rational transformation of a $d$--tuple of commuting nonselfadjoint operators.

In the first paragraph of the chapter 1 we prove classical results of the elimination theory in a framework
of Vandermonde vectors. In the rest of the chapter we describe an image of a line under a rational
transformation via a notion of Bezout matrix and give brief review of a classification of determinantal
representations of real plane algebraic curves.

In the chapter 2 we introduce the notion of Vandermonde vector on a plane algebraic curve and generalize
notions of Sylvester matrix and Bezout matrix. We formulate and prove analogues of the results of the
classical elimination theory for the elimination theory along a plane algebraic curve. We describe an image
of a plane algebraic curve given by a determinantal representation under a rational transformation via the
notion of generalized Bezout matrix.
\section{Preliminaries}
\subsection{Classical elimination theory}
Let us recall the main goal of elimination theory. Given $n+1$ (nonhomogeneous) polynomials in $n$
variables we want to find necessary and sufficient conditions (in terms of the coefficients) for these
polynomials to have a common zero (and furthermore to determine the number of common zeroes, counting
multiplicities, if they exist), see \cite{M}.

In the classical case we consider two (nonhomogeneous) polynomials in one variable. For such a pair of
polynomials one may construct so--called Sylvester matrix and Bezout matrix. Entries of these matrices
depend only on coefficients of the polynomials. Determinants of these matrices equal to zero if and only if
the polynomials have a common zero.

We will consider classical elimination theory in a framework of Vandermonde vectors. This approach allows
us to give new proofs of classical results of elimination theory for two polynomials in one variable. The
advantage of this approach is that these proofs can be generalized for elimination theory for two
polynomials in two variables along a plane algebraic curve.

Let $V_n(x)$ be the Vandermonde vector of the length $n$:
$$V_n(x) =
\left(
\begin{array}{c}
1 \\
x \\
\dots \\
x^{n-1}
\end{array}
\right) =
\left( x^i \right)_{i=0}^{n-1}$$

\begin{thm}
If $x_1, x_2, \dots , x_n$ are pairwise distinct then vectors \\
$V_n(x_1), V_n(x_2), \dots , V_n(x_n)$ are linearly independent.
\end{thm}
{\bf Proof} Let us consider $n \times n$ matrix
$$\left(V_n(x_1), V_n(x_2), \dots , V_n(x_n)\right) =
\left(
\begin{array}{cccc}
1 & 1 & \dots & 1 \\
x_1 & x_2 & \dots & x_n \\
\dots & \dots & \dots & \dots \\
x_1^{n-1} & x_2^{n-1} & \dots & x_n^{n-1}
\end{array}
\right)$$

Let us suppose that $p = (p_0, p_1, \dots , p_{n-1})$ is a row--vector from the left kernel of this matrix.
Then the polynomial $p(x) = p_0 + p_1 x + \dots +p_{n-1} x^{n-1}$ has $n$ zeroes. The degree of this
polynomial is less than $n$. Hence, the kernel of this matrix is trivial and vectors $V_n(x_1), V_n(x_2),
\dots , V_n(x_n)$ are linearly independent. Theorem is proved.

We will call $V_n(x)$ the Vandermonde vector of order zero. It is natural to define the Vandermonde vector
of higher orders. For higher orders the definition is: $V^k_n(x) =
\left(\frac{d^k}{dx^k}x^i\right)_{i=0}^{n-1}$ We will call $V^k_n(x)$ the Vandermonde vector of order $k$.

\begin{thm}
If $i_1+i_2+\dots+i_m=n-1$ and $x_1, x_2, \dots , x_m$
are pairwise distinct then vectors
$V_n(x_1), V^1_n(x_1), \dots , V^{i_1}_n(x_1),
V_n(x_2), V^1_n(x_2), \dots , V^{i_2}_n(x_2),\\
\dots, V_n(x_m), V^1_n(x_m), \dots , V^{i_m}_n(x_m)$
are linearly independent.
\end{thm}
The prove is the same as of the above theorem.

Let us consider a (nonhomogeneous) polynomial of degree $n$ in one variable:
$$p(x) = p_0 + p_1 x + p_2 x^2 \dots + p_n x^n = \sum_{i=0}^n p_i x^i =
p V_{n+1}(x)$$ It is clear that a point $x$ is a zero of the polynomial $p(x)$ if and only if the
row--vector of coefficients $p=(p_0,p_1,\dots,p_n)$ multiplied by Vandermonde vector at the point $x$ is
equal to zero: $pV_{n+1}(x)=0$.

Let us consider the $n \times 2n$ matrix
$$T(p)=
\left(
\begin{array}{cccccccc}
p_0 & p_1 &\dots& p_n & 0 &\dots&\dots& 0 \\
 0  & p_0 & p_1 &\dots& p_n & 0 &\dots& 0 \\
\dots & \dots & \dots & \dots & \dots & \dots & \dots & \dots \\
0 & \dots & 0 & p_0 & p_1 &\dots& p_n & 0 \\
0 & \dots & \dots & 0 & p_0 & p_1 & \dots & p_n
\end{array}
\right)$$

The $k$--th string of this matrix is the string of coefficients of the polynomial $x^{k-1} p(x)$. We will
call $T(p)$ the matrix of shifts of $p$.
\begin{lmm}
The kernel of the matrix of shifts of a polynomial is generated by Vandermonde vectors in zeroes of this
polynomial counting with multiplicities.
\end{lmm}
That is if the polynomial $p(x)$ has $n$ distinct zeroes $x_0, x_1,
\dots , x_n$ then \\
$\ker T(p)= span \left( V_{2n}(x_1),\dots,V_{2n}(x_n) \right)$ and if the polynomial $p(x)$ has zeroes\\
$x_0, x_1, \dots , x_m$ with multiplicities $i_0,i_1, \dots ,i_m$, $i_0+i_1+\dots+i_m=n$ then\\
$\ker T(p)= span (V_{2n}(x_1),\dots,V^{(i_1)}_{2n}(x_1),
V_{2n}(x_2),\dots,V^{(i_2)}_{2n}(x_2),\dots, \\
V_{2n}(x_m),\dots,V^{(i_m)}_{2n}(x_m))$.

Let us consider two polynomials $p$ and $q$ of degree $n$ and
$2n \times 2n$ matrix
$$S(p,q)=
\left(
\begin{array}{c}
T(p) \\
T(q)
\end{array}
\right)$$

This matrix is called Sylvester matrix of two polynomials. From the Theorem~1.2 we conclude the next
classical theorem.
\begin{thm}
The dimension of the kernel of Sylvester matrix of two polynomials is equal to the number of common zeroes
of these polynomials.
\end{thm}
\begin{cor}
Two polynomials in one variable have common zero if and only if the determinant of their Sylvester matrix
is equal to zero.
\end{cor}
The determinant of the Sylvester matrix is called a resultant.

\begin{lmm}
For every two polynomials in one variable
$p(x)$ and $q(x)$ of degree $n$
there exists uniquely determined $n \times n$ symmetric matrix
$B(p,q) = \left( b_{ij} \right)_{i,j=0}^n$ such that
$p(x)q(y) - q(x)p(y) = \sum b_{ij}x^i(x-y)y^j$.
\end{lmm}
This matrix is called Bezout matrix of two polynomials. \\
{\bf Proof} The statement of the lemma follows from the
decompositions:\\
$p(x)q(y) - q(x)p(y) = \sum p_iq_j(x^iy^j - x^jy^i)=$\\
$\sum_{i>j} p_iq_jx^j(x^{i-j} - y^{i-j})y^j -
\sum_{i<j} p_iq_jx^i(x^{j-i} - y^{j-i})y^i =$\\
$\sum_{i>j} p_iq_jx^j \left( \sum_{k=1}^{i-j} x^{k-1}(x - y)y^{i-j-k}
\right )y^j -$\\
$\sum_{i<j} p_iq_jx^i \left( \sum_{k=1}^{j-i} x^{k-1}(x - y)y^{j-i-k}
\right )y^i$
\begin{cor}
$p(x)q(y) - q(x)p(y) = V^T_n(x)(x-y)B(p,q)V_n(y)$.
\end{cor}

The following theorem is due to N. Kravitsky \cite{LMKV}, p.135.

\begin{thm} For any polynomials $p$, $q$, $f$and $g$ of degree $n$
their Sylvester matrix and Bezout matrix are connected by the next identity: \\
$S^T(p,q)
\left(
\begin{array}{cc}
0 & -B(f,g) \\
B(f,g) & 0
\end{array}
\right)
S(p,q) =\\
S^T(f,g)
\left(
\begin{array}{cc}
0 & B(p,q) \\
-B(p,q) & 0
\end{array}
\right)S(f,g)$
\end{thm}
{\bf Proof} It suffices to show that the equality is true if multiplied from the left side by $V^T_{2n}(x)$
and from the right by $V_{2n}(y)$. The left--hand side equals then $V^T_{2n}(x)S^T(p,q) \left(
\begin{array}{cc}
0 & -B(f,g) \\
B(f,g) & 0
\end{array}
\right)
S(p,q)V_{2n}(y)=$\\
$\left ( p(x)V^T_n(x), q(x)V^T_n(x) \right )
\left(
\begin{array}{cc}
0 & -B(f,g) \\
B(f,g) & 0
\end{array}
\right)
\left(
\begin{array}{c}
p(y)V_n(y) \\
q(y)V_n(y)
\end{array}
\right)=$\\
$-\left( p(x)q(y) - q(x)p(y)\right) V^T_n(x) B(f,g) V_n(y)=$\\
$-\left( p(x)q(y) - q(x)p(y)\right) \left( f(x)g(y) -
f(x)g(y)\right)(x-y)^{-1}$\\
The expression obtained preserves its value if we interchange $(p,q)$  with $(f,g)$. At the same time,
interchanging $(p,q)$ with $(f,g)$ at the left--hand side of the expression in the theorem leads us to
right--hand side. Hence, both sides are equal. Theorem is proved.

Setting $f(x)=1, g(x)=x^n$ we get the following identity.
\begin{cor}
$S^T(p,q)
\left(
\begin{array}{cc}
0 & -J_n \\
J_n & 0
\end{array}
\right)
S(p,q) =
\left(
\begin{array}{cc}
0 & B(p,q) \\
-B(p,q) & 0
\end{array}
\right)$, \\
where $J_n =
\left(
\begin{array}{cccc}
{} & {} & {} & 1 \\
{} & {} & 1 & {} \\
{} & \dots & {} & {} \\
1 & {} & {} & {}
\end{array}
\right)$.
\end{cor}
\begin{cor}
The dimension of the kernel of Bezout matrix of two polynomials $p$ and $q$ of degree $n$ is equal to the
number of common zeroes of these polynomials and $|\det S(p,q)|=|\det B(p,q)|$.
\end{cor}
The determinant of the Bezout matrix is called a bezoutian.
\subsection{Rational transformation and Bezout matrices}
The notion of the Bezout matrix allows us to describe explicitly an image of the one--dimensional
projective line ${\CP}$ to two--dimensional projective plane ${\CP}^2$ under a rational transformation.

Let us consider three homogeneous polynomials in two variables $p_0(x_0,x_1)$,\\ $p_1(x_0,x_1)$ and
$p_2(x_0,x_1)$. These polynomials define a rational transformation\\ $r: {\CP} \to {\CP}^2$ by
$$r(x_0,x_1)=(p_0(x_0,x_1), p_1(x_0,x_1), p_2(x_0,x_1))$$

Next theorem is due to Kravitsky, \cite{K}.
\begin{thm}
The image of ${\CP}$ under the rational transformation $r$ is a curve given by the next determinantal
representation:
$$\det(x_0B(p_0,p_1) + x_1B(p_2,p_0) + x_2B(p_1,p_2)) = 0$$
\end{thm}
\subsection{Determinantal representations of algebraic curves}

Let $C$ be a real projective plane curve defined by a polynomial ${\Delta}(x,y)$ of degree $m$. We will say
that $D=(D_0+xD_1+yD_2)$ is a determinantal representation of $C$ if
$${\Delta}(x,y)=\det(D_0+xD_1+yD_2)$$
Two determinantal representations $D$ and $D'$ are called Hermitian equivalent if there exists a complex $m
\times m$ matrix $P$ such that $D'=PDP^*$. We want to describe equivalence classes of determinantal
representations of $C$.

Let $D$ be a determinantal representation of $C$. For each point $(x,y)$ on $C$ consider $coker D(x,y)=\{v:
vD(x,y)=0\}$, where $v$ is $m$--dimensional row--vector. It can be shown that if $(x,y)$ is a regular point
of $C$ then $\dim coker D(x,y)=1$. Assume now $C$ is a smooth curve. It follows that $coker D$ is a line
bundle on $C$; more precisely, we define $coker D$ to be the subbundle of the trivial bundle of rank $m$
over $C$, whose fiber at the point $(x,y)$ is $coker D(x,y)$. Clearly, if two determinantal representations
$D$ and $D'$ of $C$ are equivalent, then the corresponding line bundles $coker D$ and $coker D'$ are
isomorphic. Conversely it turned out that if the line bundles corresponding to two determinantal
representations of $C$ are isomorphic, then the determinantal representations are equivalent up to sign.
The description of determinantal representations has been thus reduced to the description of certain line
bundles on $C$.

$C$ is a compact Riemann surface of genus $g$, where $g=\frac{(n-1)(n-2)}{2}$. Choosing a canonical
integral homology basis on $C$ and the corresponding normalized basis for holomorphic differentials, we
obtain the period lattice ${\Lambda}$ in ${\C}^m$. The Jacobian variety $J(C)={\C}^m / {\Lambda}$; it is a
$g$--dimensional complex torus. The Abel--Jacobi map ${\mu}$ associates to every line bundle $L$ on $C$ a
point ${\mu}(L)$ in $J(C)$. Furthermore the isomorphism class of $L$ is determined by two invariants: the
degree $\deg L$ of $L$, an integer, and the point ${\mu}(L)$ in $J(C)$.

Some important geometrical properties of the line bundle can be expressed analytically in terms of the
corresponding point in the Jacobian variety through the use of the Riemann's theta function ${\theta}(z)$.
${\theta}(z)$ is an entire function on ${\C}^g$ determined by the period lattice ${\Lambda}$. ${\theta}(z)$
is quasiperiodic with respect to ${\Lambda}$: when $z$ is translated by a vector in ${\Lambda}$,
${\theta}(z)$ is multiplied by a non--zero number, so that we can talk about the zeroes of ${\theta}(z)$ on
$J(C)$.

It can be shown that if $L=coker D$, where $D$ is a determinantal representation of $C$, then $\deg
L=-\frac{n(n-1)}{2}$. One can determine necessary and sufficient conditions on a line bundle $L$ of degree
$\frac{n(n-1)}{2}$ to be the cokernel of a determinantal representation of $C$, and translating them into
conditions on the corresponding point in the Jacobian variety yields.
\begin{thm}
$C$ possesses determinantal representation. There is a one--to--one correspondence between equivalence
classes, up to a sign, of determinantal representations $D$ of $C$ and points ${\lambda}$ of $J(C)$
satisfying ${\lambda} + \bar {\lambda} = e$ and ${\theta}({\lambda}) \neq 0$. The correspondence is given
by ${\lambda} = {\mu}(coker D(m-2))+k$.
\end{thm}
The use of the twisted line bundle $coker D(m-2)$ instead of $coker D$ and the translation of the point in
$J(C)$ by the so--called Riemann's constant $k$ are technical details. $e \in {\C}^g$ is a half--period
($2e \in {\Lambda}$) explicitly determined by the topology of the set of real points $C_{\rm {\bf R}}
\subset C$. Note that since $C$ is a real curve, the period lattice ${\Lambda}$ is invariant under complex
conjugation, and the conjugation descends to $J(C)= {\C}^m / {\Lambda}$, so the equation ${\lambda} + \bar
{\lambda} = e$ makes sense there.

This theorem is due to Vinnikov, \cite{V}.
%%%%%%%%%%%%%%%%%%%%%%%%%%%%%%%%%%%%%%%%%%%%%%%%%%%%%%%%%%%%%%%%%%%%%%%%%%%%%%%%
\section{Elimination theory on a curve}
\subsection{Vandermonde vectors on a curve}
It is an unfortunate fact that the easy, useful and beautiful constructions of the single variable case do
not generalize in any straightforward fashion to the case of two or more variables. Some generalizations of
the classical results using much more sophisticated algebraic techniques have been obtained in recent years
by Gelfand, Kapranov, Zelevinsky, see \cite{GKZ}, and Jouanolou, see \cite{J}. We shall follow a different,
somewhat ``asymmetrical'' approach. Namely, we shall choose one of our three polynomials and view it as
defining a plane algebraic curve.

Our goal becomes to determine, in terms of coefficients of two polynomials in two variables, whether these
two polynomials have a common zero on an algebraic curve defined by a third polynomial, and if they do to
determine the number of common zeroes, counting multiplicities.

To solve this problem it turns out to be essential to consider a plane algebraic curve given by a
determinantal representation
$${\Delta}(x_0,x_1,x_2)=\det(x_0D_0+x_1D_1+x_2D_2)$$
where $D_0$, $D_1$ and $D_2$ are $m \times m$ constant matrices. In this thesis we shall assume that the
(homogeneous) polynomial ${\Delta}(x_0,x_1,x_2)$ is irreducible, though many results can be generalized to
the reducible case as well. We shall also assume that ${\Delta}(x_0,x_1,x_2)$ has real coefficients and
$D_0$, $D_1$ and $D_2$ are complex hermitian matrices; it allows to dispose of certain considerations of
duality which will be necessary to handle the case when the reality assumption is dropped.

In the classical case we have considered the Vandermonde vectors in one variable: $V_n(x) = \left( x^i
\right)_{i=0}^{n-1}$ or, in homogeneous coordinates:
$$V_n(x_0, x_1) =
\left(
\begin{array}{c}
x_0^{n-1} \\
x_0^{n-2}x_1 \\
x_0^{n-3}x_1^2 \\
\dots \\
x_1^{n-1}
\end{array}
\right) =
\left(x_0^{i_0}, x_1^{i_1} \right)_{i_0+i_1=n-1} =
\left(x^i \right)_{|i|=n-1}$$

The determinantal representation of a curve provides us with an additional structure. In the case of a
plane irreducible curve for every point of the curve $x=(x_0,x_1,x_2)$ there is a nontrivial subspace
$\ker(x_0D_0 + x_1D_1 + x_2D_2)$ (one--dimensional, because of the irreducibility of
${\Delta}(x_0,x_1,x_2)$, except at the singular points). For given point $(x_0,x_1,x_2)$ we chose an
arbitrary (non--trivial) vector from the subspace $\ker(x_0D_0 + x_1D_1 + x_2D_2)$. Vandermonde vectors on
a curve (equipped with a determinantal representation) is given by
$$V_n(x,e) =
\left(
\begin{array}{c}
x_0^{n-1}e \\
x_0^{n-2}x_1e {\quad} x_0^{n-2}x_2e \\
\dots {\quad} \dots {\quad} \dots {\quad} \dots {\quad} \dots \\
x_0x_1^{n-2}e {\quad} \dots {\quad} \dots {\quad} \dots {\quad} \dots
{\quad} x_0x_2^{n-2}e \\
x_1^{n-1}e {\quad} x_1^{n-2}x_2e {\quad} \dots {\quad} \dots {\quad}
\dots {\quad} x_1x_2^{n-2}e
{\quad} x_2^{n-1}e
\end{array}
\right) =$$
$$\left( x_0^{i_0}x_1^{i_1} x_2^{i_2} e \right)_{i_0+i_1+i_2=n-1} =
\left( x^i e \right)_{|i|=n-1}$$

To define Vandermonde vectors of higher order at the point $(x_0,x_1,x_2)$ we fix a local coordinate $t$,
fix a local lifting of homogeneous coordinates from ${\CP}^2$ to ${\C}^3$~\verb@\@~$\{0\}$ and consider
$x_0(t)$, $x_1(t)$, $x_2(t)$ and $e(t) \in {\C}^m$ such that $(x_0(t)D_0 + x_1(t)D_1 + x_2(t)D_2)e(t)=0$.
In this notations $V_n(x) = \left(x_0^{i_0}(t) x_1^{i_1}(t) x_2^{i_2}(t) e(t)\right)|_{t=0} =
(x^i(t)e(t))|_{t=0}$

Now we may define Vandermonde vectors of higher order analogously to the classical case: $V_n^k(x,e)=\left(
\frac{d^k}{dt^k} x(t)e(t)\right)|_{t=0}$. This definition is correct for regular points of a curve. For
singular point of a curve we consider desingularizing Riemann surface with corresponding line bundle and
for every preimage of a point we construct in the same manner Vandermonde vectors of corresponding
multiplicities. We will call $V_n^k(x)$ the Vandermonde vector on a curve at the point $x$ of order $k$.
Vandermonde vector on a curve depend on the choice of a vector $e$ and on the choice of a local coordinate.
\begin{lmm}
For any natural number $k$ and for any point on a curve $x$ the subspace generated by Vandermonde vectors
at the point $x$ of first $k+1$ orders $V_n(x,e), V^1_n(x,e), \dots, V^k_n(x,e)$ does not depend on the
choice of a vector $e$ and on the choice of a local coordinate.
\end{lmm}
The statement of the lemma follows from the definition of Vandermonde vector of higher order and the rule
for differentiation of a product. When the choice of a vector $e$ is not essential we shall use the notion
$V_n(x)$ instead of $V_n(x,e)$.

In the classical elimination theory every Vandermonde vector $V_n(x_0,x_1)$ belongs to ${\C}^n$ and $n$
Vandermonde vectors generate the whole space ${\C}^n$.

It is clear that a Vandermonde vector $V_n(x,e)$ on a curve belongs to the space ${\bf{\rm
C}}^{m\frac{n(n+1)}{2}}$, where $m$ is the degree of ${\Delta}(x_0,x_1,x_2)$. We will denote this space by
$W_n$ and will call it the blown space. As we shall see later unlike the classical case all the Vandermonde
vectors on a curve generate not the whole blown space $W_n$ but only a certain subspace of $W_n$. Let us
consider a subspace
$$V_n = \{ (w_{j_0,j_1,j_2}) \in W_n :
D_0w_{j_0+1,j_1,j_2} + D_1w_{j_0,j_1+1,j_2} + D_2w_{j_0,j_1,j_2+1}=0 \}$$ where $j_0+j_1+j_2=n-1$. We will
call $V_n$ the principal subspace. It is clear that all Vandermonde vectors belong to $V_n$, that is
$V_n(x) \in V_n$ for every point $x=(x_0,x_1,x_2)$.
\begin{thm}
The dimension of the principal subspace $V_n$ equals to $nm$.
\end{thm}
{\bf Proof} We will prove this theorem by induction. For $n=1$ it is trivial: $V_1 = \{ w \in W_1 \}$ and
$W_1={\C}^m$.

Let us suppose now that the statement of the theorem is true for all $i < n$.
Let us consider three subspaces of $V_n$:\\
$U_1=\{ (w_{j_0,j_1,j_2}) \in W_n :
D_0w_{j_0+1,j_1,j_2} + D_1w_{j_0,j_1+1,j_2} + D_2w_{j_0,j_1,j_2+1}=0$
for $j_1 >0$,
and $w_{j_0,j_1,j_2} = 0$ for $j_1 = 0 \}$,\\
$U_2=\{ (w_{j_0,j_1,j_2}) \in W_n :
D_0w_{j_0+1,j_1,j_2} + D_1w_{j_0,j_1+1,j_2} + D_2w_{j_0,j_1,j_2+1}=0$
for $j_2 >0$,
and $w_{j_0,j_1,j_2} = 0$ for $j_2 = 0 \}$,\\
$U_3=\{ (w_{j_0,j_1,j_2}) \in W_n :
D_0w_{j_0+1,j_1,j_2} + D_1w_{j_0,j_1+1,j_2} + D_2w_{j_0,j_1,j_2+1}=0$
for $j_1=j_2=0$,
and $w_{j_0,j_1,j_2} = 0$ otherwise $\}$.

It is clear that $U_1 \cong U_2 \cong V_{n-1}$ and $U_3 \cong V_2$. Let us notice that $U_1 \cap U_2 \cap
U_3 = \emptyset$ and $U_1 \cup U_2 \cup U_3 = V_n$. Therefore $\dim V_n = \dim U_1 + \dim U_2 + \dim U_3 -
\dim(U_1 \cap U_2) - \dim(U_1 \cap U_3) - \dim(U_2 \cap U_3)$. It is clear that $(U_1 \cap U_2) \cong
V_{n-2}$ and $(U_1 \cap U_3) \cong (U_2 \cap U_3) \cong V_1$. Therefore $\dim(U_1 \cap U_2) = (n-2)m$ and
$\dim(U_2 \cap U_3) = \dim(U_1 \cap U_3) = m$. Hence, $\dim V_n = (n-1)m + (n-1)m + 2m - (n-2)m - m - m =
nm$. Theorem is proved.

Let us consider $n$ pairwise distinct linear polynomials $L_i(x)=a_ix_0 + b_ix_1 + c_ix_2$ such that
$L_i(x)$ and $L_j(x)$ do not have common zeroes on the curve for $i \neq j$.
\begin{lmm}
Vandermonde vectors in zeroes of the polynomial $\prod_{i=1}^n L_i(x)$ on a curve are linearly independent
and generate the principal subspace $V_n$.
\end{lmm}
{\bf Proof} We will prove this lemma by induction. For $n=1$ there are $m$ zeroes of $L_1(x)$ on the curve:
$x_{11},\dots,x_{1m}$, where $x_{1j}=(x_{1j,0},x_{1j,1},x_{1j,2}) \in {\CP}^2$. For every $x_{1j}$ there
exists corresponding vector $e_{1j} \in \ker(x_{1j,0}D_0+x_{1j,1}D_1+x_{1j,2}D_2)$. Let us note that
vectors $e_{11},\dots,e_{1m}$ are linearly independent as eigenvectors of a matrix pencil
$(x_0D_0+x_1D_1+x_2D_2)|_{L_1(x)=0}$, see \cite{BGR}, p14-16. Linearly independent vectors
$e_{11},\dots,e_{1m}$ generate $V_1={\C}^m$.

Let us suppose now that the statement of the lemma is true for every $i < n$.
We will consider three subspaces of $W_n$:\\
$U_1=span(V_n(x_{11}),\dots,V_n(x_{1m})$,\\
$U_2=span(V_n(x_{n1}),\dots,V_n(x_{nm})$ and\\
$U_3=span(V_n(x_{21}),\dots,V_n(x_{2m}),\dots,V_n(x_{(n-1)1}),\dots,V_n(x_{(n-1)m}))$.

By our assumption vectors $V_{n-1}(x_{11}),\dots,V_{n-1}(x_{(n-1)m})$ are linearly independent. Therefore
$\dim(U_3 \cap U_1)=0$ and, analogously, $\dim(U_3 \cap U_2)=0$. Now we may compute $\dim
span(V_n(x_{11}),\dots,V_n(x_{nm}) = \dim (U_1 + U_2 + U_3) = \dim ((U_1 + U_2) + U_3) = \dim (U_1 + U_2) +
\dim U_3 - \dim ((U_1 + U_2) \cap U_3) = \dim (U_1 + U_2) + \dim U_3 - \dim ((U_1 \cap U_3) + (U_2 \cap
U_3)) = 2m + (n-2)m = nm$. Hence, vectors $V_n(x_{11})$, $\dots,$ $V_n(x_{nm})$ are linearly independent.
It is clear that $V_n(x_{ij}) \in V_n$ for every $i$ and $j$. By previous theorem $\dim V_n = nm$. Hence,
vectors $V_n(x_{11}),\dots,V_n(x_{nm})$ generate the principal subspace $V_n$. Lemma is proved.

The blown space and the principal subspace play a key role in the elimination theory on a plane algebraic
curve. Namely, all the constructions of the classical elimination theory can be formally generalized in
terms of the blown spaces. After this formal generalization we restrict the new ``blown'' constructions to
the corresponding principal subspaces and obtain complete analogues of the classical results.
%%%%%%%%%%%%%%%%%%%%%%%%%%%%%%%%%%%%%%%%%%%%%%%%%%%%%%%%%%%%%%%%%%%%%%%%%%%%%%%%
\subsection{Generalized Sylvester and Bezout matrices}
Our goal in this paragraph is to define properly Sylvester and Bezout matrices for a pair of homogeneous
polynomials in three variables on a curve given by a determinantal representation.

As in the previous paragraph we consider the blown space $W_n$ and the principal subspace $V_n$. For a pair
of polynomials we will determine formal analogues of classical Sylvester and Bezout matrices in $W_n$ and
then will restrict these analogues on $V_n$.

Let $p(x)=p(x_0,x_1,x_2)$ be a homogeneous polynomial in three variables of degree $n$: $p(x) =
p(x_0,x_1,x_2) = \sum_{|i|=n} p_ix^i = \sum_{i_0+i_1+i_2=n} p_{i_0,i_1,i_2}x_0^{i_0}x_1^{i_1}x_2^{i_2}$.

Let us consider a $\frac{n(n+1)}{2} \times \frac{2n(2n+1)}{2}$ matrix $T(p)$, where $(i_0,i_1,i_2)$--th
string is the string of coefficients of polynomial $x^ip(x) = x_0^{i_0}x_1^{i_1}x_2^{i_2}p(x_0,x_1,x_2)$.
We will denote the matrix $T(p) \otimes I$ by $\bar T(p)$, where $I$ is the unit $m \times m$ matrix.
Multiplying $\bar T(p)$ by Vandermonde vectors we obtain main property of $\bar T(p)$.
\begin{lmm}
$\bar T(p) V_{2n}(x,e) =  p(x)V_n(x,e)$
\end{lmm}
\begin{cor}
$\bar T(p) V_{2n} \subset V_n$
\end{cor}
Let us define $T'(p) = \bar T(p)|_{V_{2n}}$ We will call the matrix $T'(p)$ the matrix of shifts of
polynomial $p$. $T'(p)$ is $nm \times 2nm$ matrix. Next lemma follows from the definition of the matrix of
shifts.
\begin{lmm}
For every polynomial $p$ of degree $n$ the rank of the matrix of shifts of this polynomial equals to $nm$:
{\,} $rk T'(p) = nm$
\end{lmm}
Let us consider two homogeneous polynomials in three variables $p(x_0,x_1,x_2)$ and $q(x_0,x_1,x_2)$ of
degree $n$ and $2nm \times 2nm$ matrix
$$S(p,q) =
\left(
\begin{array}{c}
T'(p) \\
T'(q)
\end{array}
\right)$$ We will call the matrix $S(p,q)$ the Sylvester matrix of polynomials $p$ and $q$ along the
algebraic curve given by a determinantal representation ${\Delta}$. Later we will prove that the
determinant of $S(p,q)$ equals to zero if and only if $p$ and $q$ have a common zero on a curve.
\begin{lmm}
For every two homogeneous polynomials in three variables \\
$p(x_0,x_1,x_2)$ and $q(x_0,x_1,x_2)$ of degree $n$ there exist three $\frac{n(n+1)}{2} \times
\frac{n(n+1)}{2}$ symmetric matrices ${\beta}^{10} = (b^{10}_{ij})$, ${\beta}^{20} = (b^{20}_{ij})$ and
${\beta}^{12} = (b^{12}_{ij})$ such that $p(x_0,x_1,x_2)q(y_0,y_1,y_2) - q(x_0,x_1,x_2)p(y_0,y_1,y_2) =$
$$\sum_{|i|,|j|=n}
b^{10}_{ij} x^i (x_1y_0-x_0y_1) y^j +
b^{20}_{ij} x^i (x_2y_0-x_0y_2) y^j +
b^{12}_{ij} x^i (x_1y_2-x_2y_1) y^j $$
\end{lmm}
To simplify notations we will prove nonhomogeneous analogue of this lemma.
\begin{lmm}
For every two nonhomogeneous polynomials in two variables \\
$p(x_{1},x_{2})$, $q(x_{1},x_{2})$ of degree $n$ there exist three $\frac{n(n+1)}{2} \times
\frac{n(n+1)}{2}$ symmetric matrices ${\beta}^{10} = (b^{10}_{ij})$, ${\beta}^{20} = (b^{20}_{ij})$ and
${\beta}^{12} = (b^{12}_{ij})$ such that $p(x_1,x_2)q(y_1,y_2) - q(x_1,x_2)p(y_1,y_2) =$
$$\sum_{0 \leq i_1+i_2, j_1+j_2 \leq n}
b^{10}_{ij} x^i (x_1-y_1) y^j +
b^{20}_{ij} x^i (x_2-y_2) y^j +
b^{12}_{ij} x^i (x_1y_2-x_2y_1) y^j $$
\end{lmm}
{\bf Proof } From the linearity of this decomposition it follows that it is sufficient to prove the lemma
for monomials $p(y_{1},y_{2}) = y_1^{k_1}y_2^{k_2}$, $q(y_{1},y_{2}) = y_1^{l_1}y_2^{l_2}$. Let us suppose
that $k_1 \geq l_1$, $k_2 \leq l_2$.
In this case \\
$p(x_{1},x_{2})q(y_{1},y_{2}) - q(x_{1},x_{2})p(y_{1},y_{2})=
x_1^{k_1}x_2^{k_2}y_1^{l_1}y_2^{l_2} -
x_1^{l_1}x_2^{l_2}y_1^{k_1}y_2^{k_2}=$\\
$x_1^{l_1}x_2^{k_2} (x_1^{k_1-l_1}y_2^{l_2-k_2} -
x_2^{l_2-k_2}y_1^{k_1-l_1}) y_1^{l_1}y_2^{k_2}$.\\
Let us denote $k_1-l_1=m_1$ and $l_2-k_2=m_2$. If $m_1=m_2=m$ then\\
$x_1^my_2^m - x_2^my_1^m = x_1^{m-1}(x_1y_2-x_2y_1)y_2^{m-1} + x_1^{m-2}x_2(x_1y_2-x_2y_1)y_1y_2^{m-2} +
\dots + x_1x_2^{m-2}(x_1y_2-x_2y_1)y_1^{m-2}y_2 +
x_2^{m-1}(x_1y_2-x_2y_1)y_1^{m-1}$\\
If, for example, $m_1 > m_2$ then \\
$x_1^{m_1}y_2^{m_2} - x_2^{m_2}y_1^{m_1} =
x_1^{m_1-1}(x_1y_2-x_2y_1)y_2^{m_2-1} +
x_1^{m_1-2}x_2(x_1y_2-x_2y_1)y_1y_2^{m_2-2} + \dots +
x_1^{m_1-m_2}x_2^{m_2-1}(x_1y_2-x_2y_1)y_1^{m_2-1} +$\\
$x_2^{m_2-1}(x_1y_2-x_2y_1)y_1^{m_1-1} +
x_1x_2^{m_2-2}(x_1y_2-x_2y_1)y_1^{m_1-2}y_2 + \dots +
x_1^{m_2-1}(x_1y_2-x_2y_1)y_1^{m_1-m_2}y_2^{m_2-1} -
(x_1^{m_1-m_2}x_2^{m_2}y_1^{m_2} - x_1^{m_2}y_1^{m_1-m_2}y_2^{m_2})$\\
In the expression $(x_1^{m_1-m_2}x_2^{m_2}y_1^{m_2} - x_1^{m_2}y_1^{m_1-m_2}y_2^{m_2})$ we take
$x_1^{\max(m_1-m_2,m_2)}$ and $y_1^{\max(m_1-m_2,m_2)}$ out from parenthesis and repeat the procedure until
we reduce the degree of the expression inside parenthesis to $1$. The same procedure can be done for the
case $k_1 \geq l_1$,
$k_2 \geq l_2$ because\\
$x_1^{m_1}x_2^{m_2} - y_1^{m_1}y_2^{m_2} =
x_2^{m_2}(x_1^{m_1}-y_1^{m_1}) +
(x_1^{m_1}-y_1^{m_1})y_2^{m_2} - (x_1^{m_1}y_2^{m_2} -
x_2^{m_2}y_1^{m_1})$ and
$x_1^{m_1}-y_1^{m_1} = x_1^{m_1-1}(x_1-y_1) + x_1^{m_1-2}(x_1-y_1)y_1 +
\dots +
x_1(x_1-y_1)y_1^{m_1-2} + (x_1-y_1)y_1^{m_1-1}$\\
The lemma is proved.

On the $m\frac{n(n+1)}{2}$--dimensional blown space $W_n$ let us define a $m\frac{n(n+1)}{2} \times
m\frac{n(n+1)}{2}$ matrix $B(p,q)$:
$$B(p,q) = {\beta}^{12} \otimes D_0 + {\beta}^{10} \otimes D_1 +
{\beta}^{20} \otimes D_2$$

The principal subspace ${V_n}$ is generated by Vandermonde vectors on the curve. Let us consider the
restriction of $B(p,q)$ on the principal subspace:
$$B'(p,q) = \mathcal{P}_{V_n} B(p,q) \mathcal{P}_{V_n}$$

We will call $B'(p,q)$ the Bezout matrix of polynomials $p$ and $q$ along the algebraic curve given by the
determinantal representation ${\Delta}$.

Later we will prove that the determinant of $B'(p,q)$ equals to zero if and only if $p$ and $q$ have a
common zero on a curve.
%%%%%%%%%%%%%%%%%%%%%%%%%%%%%%%%%%%%%%%%%%%%%%%%%%%%%%%%%%%%%%%%%%%%%%%%%%%%%%%%
\subsection{Main results}
In this paragraph we prove analogues of classical theorems for generalized Sylvester and Bezout matrices.

The following theorem is due to V. Vinnikov \cite{LMKV}, p 228.
\begin{thm}
If $e \in \ker(x_0D_0 + x_1D_1 - x_2D_2)$ and\\ $h \in \ker(y_0D_0 + y_1D_1 - y_2D_2)$ then
$$e^T D_0 h (x_1y_2-x_2y_1)^{-1} = e^T D_1 h (x_1y_0-x_0y_1)^{-1} = e^T
D_2 h (x_2y_0-x_0y_2)^{-1}$$
\end{thm}
We will use notation $e^T D_0 h (x_1y_2-x_2y_1)^{-1} =$ $e^T D_1 h (x_1y_0-x_0y_1)^{-1} =$\\ $e^T D_2 h
(x_2y_0-x_0y_2)^{-1} =$ $[e,h]_{x,y}$. As usually we denote indexes of summation $(i_0,i_1,i_2)$ and
$(j_0,j_1,j_2)$ by $i$ and $j$, monomials $x_0^{i_0}x_1^{i_1}x_2^{i_2}$ and $y_0^{j_0}y_1^{j_1}y_2^{j_2}$
by $x^i$ and $y^j$ and assume that $i_0+i_1+i_2=j_0+j_1+j_2=n$. Next useful lemma follows from this theorem
and from the definition of the generalized Bezout matrix.
\begin{lmm}
If
$V(x,e)$ and $V(y,h)$
are two Vandermonde vectors on a curve then
$$V^T(x,e)B'(p,q)V(y,h)=(p(x)q(y)-q(x)p(y)) [e,h]_{x,y}$$
\end{lmm}
{\bf Proof} Let us calculate $V^T(x,e)B'(p,q)V(y,h)$.\\
$V^T(x,e)B'(p,q)V(y,h)=$ (By the definition of the Bezout matrix) $=\\
V^T(x,e)({\beta}^{12} \otimes D_0 + {\beta}^{10} \otimes D_1 +
{\beta}^{20} \otimes D_2)V(y,h)=\\
\sum_{i,j}\left(
{\beta}^{12}_{i,j}x^i y^j e^T D_0 h +
{\beta}^{10}_{i,j}x^i y^j e^T D_1 h +
{\beta}^{20}_{i,j}x^i y^j e^T D_2 h\right)=\\
\sum_{i,j}\left(
{\beta}^{12}_{i,j}x^i (x_1y_2-x_2y_1) y^j e^T D_0 h
(x_1y_2-x_2y_1)^{-1}+\right.\\
{\beta}^{10}_{i,j}x^i (x_1y_0-x_0y_1) y^j e^T D_1 h
(x_1y_0-x_0y_1)^{-1} +\\
\left.
{\beta}^{20}_{i,j}x^i (x_2y_0-x_0y_2) y^j e^T D_2
h(x_2y_0-x_0y_2)^{-1})\right)=$
(By Theorem~2.2) $=\\
\sum_{i,j}\left(
{\beta}^{12}_{i,j}x^i (x_1y_2-x_2y_1) y^j+
{\beta}^{10}_{i,j}x^i (x_1y_0-x_0y_1) y^j+
{\beta}^{20}_{i,j}x^i (x_2y_0-x_0y_2) y^j\right)\cdot$\\
$[e,h]_{x,y}=(p(x)q(y)-q(x)p(y)) [e,h]_{x,y}$ Lemma is proved.
\begin{thm}
For every four polynomials in two variables $f$, $g$, $p$
and $q$ of degree $n$ the following identity holds: \\
$S'^T(p,q)
\left(
\begin{array}{cc}
0 & -B'(f,g) \\
B'(f,g) & 0
\end{array}
\right)
S'(p,q)= \\
S'^T(f,g)
\left(
\begin{array}{cc}
0 & -B'(p,q) \\
B'(p,q) & 0
\end{array}
\right)
S'(f,g)$
\end{thm}
{\bf Proof} It suffices to show that the equality is true if multiplied from the left side by
$V^T_{2n}(x,e)$ and from the right side by $V_{2n}(y,h)$.
The left--hand side equals then\\
$V^T_{2n}(x,e)S'^T(p,q)
\left(
\begin{array}{cc}
0 & -B'(f,g) \\
B'(f,g) & 0
\end{array}
\right)
S'(p,q)V_{2n}(y,h)=$ (By Lemma 2.3)\\
$\left (p(x)V^T_n(x,e), q(x)V^T_n(x,e) \right)
\left(
\begin{array}{cc}
0 & -B'(f,g) \\
B'(f,g) & 0
\end{array}
\right)
\left(
\begin{array}{c}
p(y)V_n(y,h) \\
q(y)V_n(y,h)
\end{array}
\right)=$\\
$-(p(x)q(y) - q(x)p(y)) V^T_n(x,e) B'(f,g) V_n(y,h)=$ (By Lemma 2.7)
$=\\
-(p(x)q(y) - q(x)p(y)) (f(x)g(y) - g(x)f(y)) [e,h]_{x,y}$\\
The expression obtained preserves its value if we interchange $(p,q)$ with $(f,g)$. At the same time,
interchanging $(p,q)$ with $(f,g)$ at the left--hand side of the expression in the theorem leads us to
right--hand side. Hence, both sides are equal.
Theorem is proved.\\
Now we may formulate preliminary result.
\begin{lmm}
Let us suppose that Vandermonde vectors in zeroes of polynomial $q(x_0,x_1,x_2)$ on a curve (counting with
multiplicities) are linearly independent. Then\\
1. The kernel of the matrix of shifts of this polynomial is generated by Vandermonde vectors in zeroes of
this polynomial on a curve.\\
2. For every polynomial $p(x_0,x_1,x_2)$ the dimension of the kernel of Sylvester matrix of $p$ and $q$
along the curve is equal to the number of common zeroes of these polynomials on the curve.\\
3. For every polynomial $p(x_0,x_1,x_2)$ the dimension of the kernel of Bezout matrix of $p$ and $q$ along
the curve is equal to the number of common zeroes of these polynomials on the curve.
\end{lmm}
{\bf Proof} It follows from lemma 2.4 that the dimension of the kernel of matrix of shifts equals to $nm$.
It follows from lemma 2.4 that all $nm$ independent Vandermonde vectors in zeroes of polynomial
$q(x_0,x_1,x_2)$ belong to this kernel. Hence, the first statement of this lemma is true. To obtain the
second statement of this lemma we consider a vector $U$ from the kernel of the Sylvester matrix of $p$ and
$q$. From the definition of the Sylvester matrix it follows that this vector belongs to the kernel of the
matrix of shifts of polynomial $q$. Hence, $U=\sum a_iV_i(y_i)$, where $y_i$ are zeroes of polynomial $q$
and $V_i(y_i)$ are Vandermonde vectors in zeroes of polynomial $q$. Analogously, vector $U$ belongs to the
kernel of the matrix of shifts of polynomial $p$. That is $T'(p)U=0$. But by lemma 2.3 $T'(p)U=\sum
p(y_i)a_iV_i(y_i)$ and by assumption all $V_i(y_i)$ are linearly independent. Hence, $p(y_i)=0$ and the
second statement of this lemma is true. To prove third statement of this lemma let us assume that the curve
does not contain points $(0,1,0)$ and $(0,0,1)$. Let us consider polynomials $x_1^n$ and $x_2^n$. From the
definition of the generalized Bezout matrix it follows that $B'(x_1^n, x_2^n)$ is nondegenerate. Third
statement follows now from the second one and from theorem 2.3. If the curve contains points $(0,1,0)$ and
$(0,0,1)$ we make an invertible linear change of variables such that new curve does not contain points
$(0,1,0)$ and $(0,0,1)$, take polynomials $x_1^n$ and $x_2^n$ on the new curve and then make an inverse
change of variables. The generalized Bezout matrix of the new pair of polynomials will be nondegenerate and
we again may apply the second statement of this theorem and theorem 2.3. Lemma is proved.

It is clear now that if we will manage to prove that Vandermonde vectors in zeroes of an arbitrary
polynomial are linearly independent then we are done.
\begin{lmm}
For every $k$ points from ${\C}$: $a_1, a_2, \dots, a_k$ there exist two numbers $c_1, c_2 \in {\C}$ such
that $c_1a_i+c_2a_j \not= 0$ for every $1 \leq i,j \leq k$.
\end{lmm}
{\bf Proof} All pairs $c_1, c_2$ that does not satisfy the statement of this lemma lies on the finite
number of lines in ${\C}^2$. We chose a pair that does not lie on these lines. Lemma is proved.
\begin{thm}
For every polynomial Vandermonde vectors in zeroes of this polynomial on a curve (counting with
multiplicities) are linearly independent.
\end{thm}
That is if
$x^0=(x_0^0,x_1^0,x_2^0),x^1=(x_0^1,x_1^1,x_2^1),\dots,x^k=(x_0^k,x_1^k,x_2^k)$
are zeroes of polynomial $p(x)$ on a curve with
multiplicities $i_0,i_1,\dots,i_k$, $i_0+i_1+\dots+i_k = nm$ then
$V_{n}(x^0), \dots, V^{i_0}_{n}(x^0),
V_{n}(x^1), \dots, V^{i_1}_{n}(x^1), \dots,
V_{n}(x^k), \dots, V^{i_k}_{n}(x^k)$ are linearly independent.\\
${\qquad}$\\
{\bf Proof} Let us consider polynomial $p(x_0,x_1,x_2)$. For simplicity we will suppose that
$p(x_0,x_1,x_2)$ has no common zeroes on the curve with the polynomial $x_0$ and that values of $x_0$ and
$x_1$ at zeroes of $p$ on the curve are pairwise distinct (by previous lemma it may be done by linear
change of variables). Let us consider $n-1$ linear polynomials $L_1(x),\dots,L_{n-1}(x)$ such that $L_i(x)$
has no common zeroes on the curve with the polynomials $x_0, x_1, p(x)$ and $L_j(x)$ for every $i \neq j$.
Let us denote $\prod_{i=1}^{n-1}L_i(x)$ by $q(x)$. Let us consider a matrix
$A=x_1B'(x_0q(x),p)+x_0B'(x_1q(x),p)$. It follows from lemma 2.2 that Vandermonde vectors in zeroes of the
polynomial $x_0q(x)$ on the curve are linearly independent. Hence, from lemma 2.8 we obtain that $\det
B'(x_0q(x),p)\not=0$ and therefore the determinant of $A$ is not identically zero. Now, as in the
definition of Vandermonde vectors of higher order, we fix a local coordinate $t$, fix a local lifting of
homogeneous coordinates from ${\CP}^2$ to ${\C}^3$~\verb@\@~$\{0\}$ and consider $x_0(t)$, $x_1(t)$ and
$x_2(t)$. In this notations $A(t)=x_1(t)B'(x_0^n(t),p(x(t)))+x_0(t)B'(x_0^{n-1}(t)x_1(t),p(x(t)))$ From the
definition of Bezout matrix it follows that Vandermonde vectors in zeroes of polynomial $p$ on a curve are
a set of null chains of $A(t)$ Hence, they are linearly independent. (For the definition and properties of
null chains, see \cite{BGR}, p.14--17). Theorem is proved.

Now, when we have proved that Vandermonde vectors in zeroes of every polynomial are linearly independent,
we may deduce from lemma 2.8 main results.
\begin{lmm}
For every polynomial the kernel of the matrix of shifts of this polynomial is generated by Vandermonde
vectors in zeroes of this polynomial on a curve.\\
$\ker T'(p) =
<V_{n}(x_0), \dots, V^{i_0}_{n}(x_0),
V_{n}(x_1), \dots, V^{i_1}_{n}(x_1), \dots,\\
V_{n}(x_k), \dots, V^{i_k}_{n}(x_k)>$, where $x_0,x_1,\dots,x_k$ are zeroes of polynomial $p(x)$ with
multiplicities $i_0,i_1,\dots,i_k$, $i_0+i_1+\dots+i_k = nm$.
\end{lmm}
\begin{thm}
The dimension of the kernel of Sylvester matrix of two polynomials along the curve is equal to the number
of common zeroes of these polynomials on the curve.
\end{thm}
\begin{cor}
Two polynomials have common zero on a curve if and only if the determinant of their Sylvester matrix along
this curve is equal to zero.
\end{cor}
\begin{thm}
The dimension of the kernel of Bezout matrix of two polynomials along the curve is equal to the number of
common zeroes of these polynomials along the curve.
\end{thm}
\begin{cor}
Two polynomials in two variables have common zero on a curve if and only if the determinant of their Bezout
matrix along this curve is equal to zero.
\end{cor}
\subsection{Rational transformation of a curve}
As was already mentioned in the paragraph~1.2 of the chapter~1 an image of the strait line under a rational
transformation can be described via the notion of the Bezout matrix. In this paragraph we consider a
rational transformation of a curve given by a determinantal representation and describe an image of a curve
via the notion of the generalized Bezout matrix.

As before, we start from a real plane algebraic curve $C$ defined by a homogeneous polynomial in three
variables ${\Delta}(x_0,x_1,x_2)$:
$$C = \{ (x_0,x_1,x_2) \in {\CP}^2: {\Delta}(x_0,x_1,x_2)=0 \}$$
and a determinantal representation of this curve:
$${\Delta}(x_0,x_1,x_2)= \det (x_0D_0 + x_1D_1 + x_2D_2)$$

Let us consider three homogeneous polynomials in three variables $p_0$, $p_1$, $p_2$ of degree $n$ and an
image $r(C)$ of the curve $C$ under the rational transformation $r=(p_0,p_1,p_2)$:
$$r(C) = \{ (p_0(x),p_1(x),p_2(x)) \in {\CP}^2: x \in C\}$$

To find a determinantal representation of the curve $r(C)$ we consider three generalized Bezout matrices:
$B'(p_0,p_1)$, $B'(p_0,p_2)$ and $B'(p_1,p_2)$.
\begin{lmm}
The curve given by determinantal representation\\
$\det (x_0B'(p_1,p_2) + x_1B'(p_2,p_0) + x_2B'(p_0,p_1))$ contains
$r(C)$.
\end{lmm}
{\bf Proof} Let $V(x,e)$ and $V(y,h)$ be two Vandermonde vectors on a curve.\\
$V^T(x,e)(p_0(x)B'(p_1,p_2) + p_1(x)B'(p_0,p_2) +
p_2(x)B'(p_0,p_1))V(y,h)=$\\
$p_0(x)V^T(x,e)B'(p_1,p_2)V(y,h) + p_1(x)V^T(x,e)B'(p_2,p_0)V(y,h) +\\
p_2(x)V^T(x,e)B'(p_0,p_1)V(y,h)=$
By Lemma 2.7 $=\\
p_0(x)(p_1(y)p_2(x) - p_2(y)p_1(x))[e,h]_{x,y} +
p_1(x)(p_2(y)p_0(x) - p_0(y)p_2(x))[e,h]_{x,y} +
p_2(x)(p_1(y)p_0(x) - p_0(y)p_1(x))[e,h]_{x,y} = 0$\\
This identity is hold for arbitrary $y=(y_0,y_1,y_2) \in C$. Hence, by
lemma 2.2\\
$(p_0(x)B'(p_1,p_2) + p_1(x)B'(p_2,p_0) + p_2(x)B'(p_0,p_1))V(x)=0$ and therefore\\
$\det (p_0(x)B'(p_1,p_2) + p_1(x)B'(p_2,p_0) + p_2(x)B'(p_0,p_1))=0$. Lemma is proved.

But if a basepoint $y$ of the rational transformation $r=(p_0,p_1,p_2)$ belongs to the curve $C$, that is
if there exists a common zero of polynomials $p_0$, $p_1$ and $p_2$ on the curve, then $\det
(x_0B'(p_1,p_2) + x_1B'(p_2,p_0) + x_2B'(p_0,p_1))=0$ because $B'(p_1,p_2)V(y)=B'(p_2,p_0)=B'(p_0,p_1)=0$.

Therefore we consider a subspace $\bar V_n$ of the principal subspace $V_n$ which is perpendicular to the
Vandermonde vectors in basepoints and restrict generalized Bezout matrices on this
subspace:\\
$\bar V_n = <V_n(y): p_0(y)=p_1(y)=p_2(y)=0>^{\perp}$, $\bar B(p_i,p_j) = \mathcal{P}_{\bar V_n}
B'(p_i,p_j) \mathcal{P}_{\bar V_n}$.

From the corollary~2.3 we conclude next theorem
\begin{thm}
The curve $r(C)$ has the determinantal representation\\
$\det (x_0 \bar B(p_1,p_2) + x_1 \bar B(p_2,p_0) +
x_2 \bar B(p_0,p_1))$.
\end{thm}
\begin{cor}
If the curve $C$ has two equivalent determinantal representations then two determinantal representations of
the curve $r(C)$ will be equivalent.
\end{cor}
The proof follows from the previous theorem and from the definition of generalized Bezout matrix.

Let us consider now two rational transformations $r=(p_0,p_1,p_2)$ and $s=(q_0,q_1,q_2)$. Here
$p_0,p_1,p_2$ are homogeneous polynomials in three variables of degree $n$ and $q_0,q_1,q_2$ are
homogeneous polynomials in three variables of degree $k$. It is clear that previous theorem allows us to
construct a determinantal representations of the curve $s(r(C))$ in two different ways:

-- we may construct the determinantal representation of the curve $C$
under the rational transformations $r \circ s$:\\
$\det (x_0 \bar B(q_1(p),q_2(p)) + x_1 \bar B(q_2(p),q_0(p)) + x_2 \bar
B(q_0(p),q_1(p)))$

-- we may construct the determinantal representation of
the curve $r(C)$ under the rational transformations $s$:\\
$\det (x_0 \bar B(p_1,p_2) + x_1 \bar B(p_2,p_0) + x_2 \bar
B(p_0,p_1))$

More precisely, we consider three homogeneous polynomials in three variables $q_0,q_1,q_2$ of degree $k$ on
the curve $r(C)$ given by the determinantal representations $\det (x_0 \bar B(p_1,p_2) + x_1 \bar
B(p_2,p_0) + x_2 \bar B(p_0,p_1))$. For this case we denote a blown space $\bar V_n^{\frac{k(k+1)}{2}}$ by
$W_k^{r(C)}$ and a principal subspace $\{ (w_{j_0,j_1,j_2}) \in W_k^{r(C)} : \bar
B(p_1,p_2)w_{j_0+1,j_1,j_2} + \bar B(p_2,p_0)w_{j_0,j_1+1,j_2} + \bar B(p_0,p_1)w_{j_0,j_1,j_2+1}=0 \}$ by
$V_k^{r(C)}$ (here $j_0+j_1+j_2=k-1$). The generalized Bezout matrices of the polynomials $q_i$ and $q_j$
along the curve $r(C)$ are denoted by $B'_{r(C)}(q_i,q_j)$ and after the restriction of these matrices on
the subspace that perpendicular to Vandermonde vectors in common zeroes of $q_0$, $q_1$ and $q_2$ on the
curve $r(C)$ we obtain the matrices $\bar B_{r(C)}(q_i,q_j)$ and the determinantal representation of the
curve $s(r(C))$: $\det (x_0 \bar B_{r(C)}(p_1,p_2) + x_1 \bar B_{r(C)}(p_2,p_0) + x_2 \bar
B_{r(C)}(p_0,p_1))$.

\begin{thm}
Two determinantal representations\\
$\det (x_0 \bar B_C(q_1(p),q_2(p)) +
x_1 \bar B_C(q_2(p),q_0(p)) +
x_2 \bar B_C(q_0(p),q_1(p)))$ and\\
$\det (x_0 \bar B_{r(C)}(p_1,p_2) + x_1 \bar B_{r(C)}(p_2,p_0) + x_2
\bar B_{r(C)}(p_0,p_1))$
are equivalent.
\end{thm}
{\bf Proof} Let us consider a map $M_i : W_{nk} \to V_n$,
$M_i(w_j) = (v_l)$ defined by a formula
$v_{l_0,l_1,l_2} = w_{i_0+l_0,i_1+l_1,i_2+l_2}$.
Here
$i=(i_0,i_1,i_2)$, $i_0+i_1+i_2=nk-n$, $j=(j_0,j_1,j_2)$,
$j_0+j_1+j_2=nk$,
$l=(l_0,l_1,l_2)$, $l_0+l_1+l_2=n$.
For every polynomial ${\alpha}(x)=\sum{\alpha}_ix^i$ of degree $nk-n$
there exists a map $S({\alpha}): W_{nk} \to V_n$ defined by a formula
$S({\alpha})=\sum{\alpha}_i M_i$. Let us note that if
${\alpha}(x)=p_0^{i_0}p_1^{i_1}p_2^{i_2}$
then $S({\alpha})$ maps $W_{nk}$ in $\bar V_n$.

Let us consider now space
$W^{r(C)}_k = \bar V_n^{\frac{k(k+1)}{2}} = \oplus \bar V_{n,i} =
\oplus \bar V_{n;i_0,i_1,i_2}$
and a map ${\tau}: W_{nk} \to W^{r(C)}_k$ defined by formula ${\tau} =
({\tau}_i) = ({\tau}_{i_0,
i_1,i_2})$,
where ${\tau}_{i_0,i_1,i_2}: W_{nk} \to \bar V_{n;i_0,i_1,i_2}$ and
${\tau}_{i_0,i_1,i_2}=S(p_0^{i_0}p_1^{i_1}p_2^{i_2})$.

From the definition of the generalized Bezout matrices it follows that\\
$B(q_i(p),q_j(p))={\tau}B_{r(C)}(q_i,q_j){\tau}^*$ for $i,j=0,1,2$.

Let us consider an action of the projection $\mathcal{P}_{\bar V_{nk}} {\tau}$ on Vandermonde vectors that
form a basis of $\bar V_{nk}$. It is clear that the images of these vectors form a basis of $\bar
V^{r(C)}_k$. Hence, $\mathcal{P}_{\bar V_{nk}} {\tau}$ is an isomorphism between $\bar V_{nk}$ and $\bar
V^{r(C)}_k$. Theorem is proved.

\end{document}